\magnification=\magstep1
\input amstex
\UseAMSsymbols
\input pictex
\vsize=23truecm
\NoBlackBoxes
\pageno=1

   \font\rmk=cmr8    \font\itk=cmti8  \font\ttk=cmtt8

   \font\gross=cmbx10 scaled\magstep1 

\def\soc{\operatorname{soc}}

\def\Hom{\operatorname{Hom}}
\def\End{\operatorname{End}}

\def\Ker{\operatorname{Ker}}

\def\bdim{\operatorname{\bold{dim}}}

\def\arr#1#2{\arrow <1.5mm> [0.25,0.75] from #1 to #2}

\vglue1truecm

\centerline{\gross The eigenvector variety of a matrix pencil.}
		          \bigskip
\centerline{Claus Michael Ringel}
	\bigskip\medskip
{\narrower\narrower Abstract.
    Let $k$ be a field and $n,a,b$ natural numbers. A matrix pencil $P$ is given by
    $n$ matrices of the same size with coefficients in $k$, say by $(b\times a)$-matrices,
    or, equivalently, by $n$ linear transformations $\alpha_i\:k^a \to k^b$ with $i=1,\dots,n$.
  We say that $P$ is {\it reduced} provided the intersection of the kernels
  of the linear transformations $\alpha_i$ is zero.
  If $P$ is a reduced matrix pencil, a vector $v\in k^a$ will be
    called an {\it eigenvector} of $P$ provided the subspace
  $\langle \alpha_1(v),\dots,\alpha_n(v) \rangle$
 of $k^b$ generated by the elements $\alpha_1(v),\dots,\alpha_n(v)$ is $1$-dimensional.
 Eigenvectors are called {\it equivalent} provided they are scalar
 multiples of each other. The set $\epsilon(P)$ of equivalence
 classes of eigenvectors of $P$ is a Zariski closed subset of the projective
 space $\Bbb P(k^a)$, thus a projective variety.
 We call it the {\it eigenvector variety} of $P$.
 The aim of this note is to show that
 any projective variety arises as an eigenvector variety of some reduced matrix pencil.\par}
								    	       	               \bigskip \medskip
													  \plainfootnote{} {\rmk 2010 \itk Mathematics Subject Classification. \rmk
		15A22,  
		14J10, 
		16G20. 
\itk Keywords: \rmk
  Matrix pencils. Eigenvectors, eigenvalues. Projective varieties.
  Kronecker modules. Quiver Grassmannians.
}

{\bf 1\. Introduction.}
     \medskip
     Let $k$ be a field and $n,a,b$ natural numbers. A {\it matrix pencil} $P$
     is given by
     $n$ matrices of the same size with coefficients in $k$, say $(b\times a)$-matrices
     or, equivalently, by $n$
     linear transformations $\alpha_i\:k^a \to k^b$ with $i=1,\dots,n$, thus we write
     $P = (a,b;\alpha_1,\dots,\alpha_n)$ and we call $(a,b)$ the
     {\it dimension vector} of $P$.
     We say that $P$ is {\it reduced} provided the intersection of the kernels
     of the linear transformations $\alpha_i$ is zero.

Let us assume that $P = (a,b;\alpha_1,\dots,\alpha_n)$ is a reduced
matrix pencil. A vector $v\in k^a$ will be
called an {\it eigenvector} of $P$ provided the subspace of $k^b$ generated by
the elements $\alpha_1(v),\dots,\alpha_n(v)$ is $1$-dimensional.
Of course, if $v$ is an eigenvector
of $P$, then any non-zero multiple of $v$ is also an eigenvector;
eigenvectors are called
{\it equivalent} provided they are multiples of each other.
The set $\epsilon(P)$ of equivalence
classes of the eigenvectors of $P$ is a Zariski closed subset of the projective
space $\Bbb P(k^a)$, thus a projective variety. We call it the {\it eigenvector variety} of $P$.
If $n = 2$ and $\alpha_1$ is the identity $(a\times a)$-matrix, then the eigenvectors of
the matrix pencil $(a,a;\alpha_1,\alpha_2)$ are just the eigenvectors of $\alpha_2$ as usually considered.

The aim of this note is to show the following.
    \medskip

{\bf Theorem 1.} {\it Let $k$ be algebraically closed.
Any projective variety arises as the eigenvector variety $\epsilon(P)$ of a reduced matrix
pencil $P$.}
       \medskip
We may require, in addition, that the reduced matrix pencil $P$ is a matrix pencil of
square matrices, see 4.2.

	\medskip
In order to prove Theorem 1, we will reformulate
the setting in terms of rings and modules, namely as
dealing with bristle varieties of Kronecker modules. But before we do this,
let us mention that the paper includes an appendix which provides further
information on the existence of eigenvectors of reduced matrix pencils.
       		       \medskip
Here is the reformulation of Theorem 1 in terms of quivers and their representations
(see, for example, [9]), thus in terms of rings and modules.
Let $K(n)$ be the $n$-Kronecker quiver; it consists of two vertices,
denoted by $1$ and $2$, and $n$ arrows
$\alpha_i\:1 \to 2$, with $1\le i \le n$.
$$
 \hbox{\beginpicture
  \setcoordinatesystem units <3cm,2.6cm>
 \put{$1$} at 0 0
 \put{$2$} at 1 0
 \arr{0.2 0.1}{0.8 0.1}
 \arr{0.2 -.1}{0.8 -.1}
 \put{$\alpha_1$} at 0.5 0.2
 \put{$\alpha_{n}$} at 0.5 -.2
 \put{$\vdots$} at  0.5 0.035
 \endpicture}
$$
The $n$-Kronecker modules are the representations of $K(n)$ over $k$.
Kronecker modules are written in the form $M =
(M_1,M_2;\alpha_1,\dots,\alpha_{n})$
where $M_1,M_2$ are vector spaces and $\alpha_i\:M_1\to M_2$ are
linear transformations, for $1\le i \le n$; the dimension vector of $M$ is
the pair $\bdim M = (\dim M_1,\dim M_2).$ If 
$(a,b;\alpha_1,\dots,\alpha_{n})$ is matrix pencil, we call $M(P) = 
(k^a,k^b;\alpha_1,\dots,\alpha_{n})$
the corresponding Kronecker module. Conversely, given a Kronecker module
$M = (M_1,M_2;\alpha_1,\dots,\alpha_{n})$, we may choose bases $\Cal B_1, \Cal B_2$
of $M_1, M_2$, respectively, so that we can write the linear transformations 
$\alpha_i$ as matrices. We obtain in this way a
corresponding matrix pencil $P(M;\Cal B_1,\Cal B_2).$ 

The $n$-Kronecker algebra $\Lambda(n) = kK(n)$ is the path algebra of $K(n)$, thus the
$\Lambda(n)$-modules are just $n$-Kronecker modules. There are precisely two
(isomorphism classes of) simple $n$-Kronecker modules, $S(1)$ and $S(2)$  
(in general, if $x$ is a vertex of an acyclic quiver, we denote by $S(x)$ the corresponding
simple module; it is defined by $S(x)_x = k$ and $S(x)_y = 0$ for all vertices $y\neq x$).
The module $S(1)$ is injective, the module $S(2)$
is projective. A $\Lambda(n)$-module $M$ is called {\it reduced} provided $S(1)$ 
is not a submodule of $M$, thus iff the corresponding matrix pencils are reduced.

In general, given a ring $R$, a {\it bristle} is by definition an
indecomposable $R$-module of length 2.
An $n$-Kronecker module is a bristle if and only if it is indecomposable with
dimension vector $(1,1)$. Such a bristle
is of the form $B(\lambda) = (k,k;\lambda_1,\dots,\lambda_{n})$, where
$(\lambda_1,\dots,\lambda_{n})$ is a non-zero element of $k^n$.
				     	      \medskip
Given a representation $M$ of a quiver and $\bold d$ a dimension vector,
one denotes by $\Bbb G_{\bold d}(M)$ the set of submodules of $M$ with dimension vector $\bold d$
and calls it a {\it quiver Grassmannian;} this is a projective variety.

Let $M$ be a reduced Kronecker module.
Since we assume that $M$ is reduced, any submodule of $M$ with dimension vector $(1,1)$ is
indecomposable, thus a bristle.
 We write $\beta(M) = \Bbb G_{(1,1)}(M)$
 for the set of bristle submodules of $M$ and call it the {\it bristle variety} of
 $M$. Of course, $\beta(M)$ is nothing else than the eigenvector variety of the
 matrix pencils corresponding to $M$:
 $$
  \beta(M) = \epsilon(P(M,\Cal B_1,\Cal B_2))
  $$
  for bases $\Cal B_1, \Cal B_2$ of $M_1, M_2$, respectively.
      \medskip
Theorem 1 may be formulated as follows:
      	      \medskip
{\bf Theorem 1$'$.} {\it Let $k$ be algebraically closed.
Any projective variety arises as the bristle variety
$\beta(M)$ of a reduced Kronecker module $M$.}
	         	  \medskip
This formulation shows that here we are in the realm of a highly contested topic:
the realization of any projective variety as the moduli space of serial modules, as studied by
Huisgen-Zimmermann and several coauthors [7,1,2],
as the moduli space of thin modules, studied by Hille [5],
and as a quiver Grassmannian, studied by Reineke [8]; see [10] and the discussions
in various internet blogs quoted there.
The main result of the present note was first
presented at the conference in honor of Jerzy Weyman's 60th birthday, April 2015, at
the University of Connecticut and then on various other occasions.
Our construction is obviously based on
the previous accounts, but the observation that it is sufficient
to consider length 2 submodules seems to have been new at that time.
For further discussions, we want to refer to Hille [6].
    	          \medskip

We say that a module is {\it bristled} provided
it is generated by modules of length $2$.
The paper [11] has been devoted to the study of bristled Kronecker modules.
There, it has been
shown that for $n\ge 3$, there is an abundance of bristled Kronecker modules
(see the appendix of the present paper).
Of course, also the main result of the present note deals with bristled Kronecker modules.
A module $M$ will be called {\it fully bristled} provided it is bristled and any simple submodule
of $M$ is contained in a bristle submodule of $M$.

There is an interesting fully bristled
$n$-Kronecker module $C$ which will be exhibited in the next section, we call it the canonical bristled
$n$-Kronecker module. The proof of Theorem 1$'$
will be established by looking at submodules of $C$.
         \bigskip\medskip
{\bf 2\. The canonical bristled module $C$.}
     \medskip
     The aim of this section is to construct a special $n$-Kronecker module $C$ which we will
     call the canonical bristled module. Here we allow that $k$ is an arbitrary field.

We define $C = (C_1,C_2;\alpha_1,\dots,\alpha_n)$ as follows:
its dimension vector is $\bdim C = (\binom{n+1}2,n)$, the space
$C_2$ has the basis elements $c_i,$ with $1\le i \le n$,
the space $C_1$ has the basis elements $c_{ij} = c_{ji}$, where $1\le i \le j \le n$,
and $\alpha_i(c_{ij}) = c_j,$ whereas $\alpha_r(c_{ij}) = 0$ in case $r\notin\{i,j\}$.
As we will see, $C$ is always an indecomposable bristled module;
we call it the {\it canonical bristled module}.

For $n=1$, the module $C$ is the only indecomposable module which is not simple;
it is just the unique bristle. For $n=2$, the canonical bristled module is the unique
indecomposable $2$-Kronecker module with dimension vector $(3,2)$ (it is the Auslander-Reiten
translate of the simple injective module $S(1)$).
Here are illustrations for $n=2$ and $n=3$:
$$
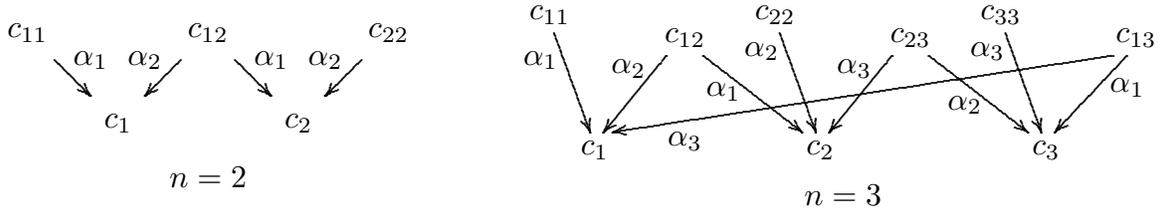
\hbox{\beginpicture
\setcoordinatesystem units <1cm,1cm>
 \put{\beginpicture
 \put{$c_{11}$} at 0 1
 \put{$c_{12}$} at 2 1
 \put{$c_{22}$} at 4 1
 \put{$c_{1}$} at 1 0
 \put{$c_{2}$} at 3 0

\arr{0.3 0.7}{0.7 0.3}
\arr{1.7 0.7}{1.3 0.3}
\arr{2.3 0.7}{2.7 0.3}
\arr{3.7 0.7}{3.3 0.3}
\put{$\alpha_1$} at 0.7 0.7
\put{$\alpha_2$} at 1.3 0.7
\put{$\alpha_1$} at 2.7 0.7
\put{$\alpha_2$} at 3.3 0.7
\put{$n=2$} at 2 -.6
\endpicture} at 0 0
\put{\beginpicture
\setcoordinatesystem units <.5cm,.5cm>
\put{$c_{11}$} at 0 3
\put{$c_{1}$} at 1 0
\put{$c_{12}$} at 3 2.5
\put{$c_{22}$} at 5 3
\put{$c_{2}$} at 6 0
\put{$c_{23}$} at 8 2.5
\put{$c_{33}$} at 10 3
\put{$c_{3}$} at 11 0
\put{$c_{13}$} at 13 2.5

\arr{0.1 2.6}{0.9 0.4}
\arr{5.1 2.6}{5.9 0.4}
\arr{10.1 2.6}{10.9 0.4}
\arr{2.6 2.1}{1.2 0.4}
\arr{7.6 2.1}{6.2 0.4}
\arr{3.4 2.1}{5.6 0.4}
\arr{8.4 2.1}{10.6 0.4}

\arr{12.8 2.1}{11.2 0.4}
\arr{12.5 2.1}{1.4 0.4}

\put{$\alpha_1$} at -.2 2
\put{$\alpha_2$} at 4.7 2.2
\put{$\alpha_3$} at 9.7 2.2
\put{$\alpha_2$} at 1.75 1.8
\put{$\alpha_3$} at 6.75 1.8
\put{$\alpha_1$} at 12.8 1.4
\put{$\alpha_1$} at 3.85 1.3
\put{$\alpha_2$} at 9.2 1
\put{$\alpha_3$} at 3 0.2
\put{$n=3$} at 6.5   -1
\endpicture} at 7 0
\endpicture}
$$
	\medskip
	{\bf (2.1)} {\it For every  non-zero
	element  $c = \sum_i \lambda_ic_i$ of $C_2$, there is a unique bristle
	$B$ with $c\in B$, it is generated by
	$d(c) = \sum_{i\le j}\lambda_i\lambda_j c_{ij}\in C_1$
	and has type $(\lambda_1:\lambda_2:\cdots:\lambda_n)$. It follows that
	$C$ is a fully bristled module, that $\dim\Hom(B,C) = 1$ for any
	bristle $B$ and that $\End(C) = k.$}
		\medskip
Proof.
Given a linear combination $d = \sum_{i\le j}\lambda_{ij}c_{ij}$, we have $\alpha_r(d) =
\sum_{i}\lambda_{ir}c_i$ where we write $\lambda_{ij} = \lambda_{ji}$ for all $i,j$.
It follows that we have $\alpha_r(d) = 0$ if and only if
$\lambda_{ir} = 0$ for $1\le i \le r$. As a consequence, $\alpha_r(d) = 0$ for all
$1\le r \le n$ implies that $d = 0.$ This shows that $C$ has no direct summand of the form
$S(1).$

Next, consider a non-zero element  $c = \sum_i \lambda_ic_i$ of $C_2$ and let
$d = d(\lambda)
= \sum_{i\le j}\lambda_i\lambda_j c_{ij}\in C_1.$ We see that $\alpha_r(d) =
\sum_{i}\lambda_{i}\lambda_r c_i = \lambda_r c$ is a multiple of $c$, and for at least
one index $r$, the multiple $\lambda_rc$ is non-zero.
Thus, the submodule generated by $d$ is the vector space with basis $c,d$. It is a bristle
and $\alpha_r(d) = \lambda_r\cdot c$ shows that this bristle is of type
$(\lambda_1:\lambda_2:\cdots:\lambda_n)$.

On the other hand, consider a submodule $U$ of $C$ which is a bristle
and assume that $c = \sum_i \lambda_ic_i$ generates its socle.
Assume that $d = \sum_{i\le j}\lambda_{ij}c_{ij}$ generates $U$.
Since $U$ is a bristle, $\alpha_j(d)$ is a multiple of $c$, for any $j$,
say $\alpha_j(d) = \mu_j\cdot c$
for some $\mu_j.$ This means that
$$
  \sum_i \lambda_{ij}c_i = \alpha_j(d) = \mu_j\cdot c = \sum_i \mu_j\lambda_i c_i,
  $$
  and therefore $\lambda_{ij} = \lambda_i\mu_j$
  for all $i,j$. Since $\lambda_{ij} = \lambda_{ji}$, we have
  $\lambda_i\mu_j = \lambda_j\mu_i$
  for all $i,j$. Now $0 \neq c = \sum_i \lambda_{i}c_i$, thus $\lambda_t\neq 0$
  for some $t$. The equality $\lambda_t\mu_j = \lambda_j\mu_t$ implies that
  $$
   \mu_j = \frac{\mu_t}{\lambda_t}\lambda_j,
   $$
   for all $j$. As a consequence,
   $$
    \lambda_{ij} = \lambda_i\mu_j = \frac{\mu_t}{\lambda_t}\lambda_i\lambda_j,
    $$
    thus $d$ is a multiple of $d(\lambda)$. This shows that the only submodule
    $U$ which is a bristle with socle generated by $c = \sum_i \lambda_ic_i$ is
    the bristle generated by $d(c)$. For later reference, let us note that
    we have shown in this way that $C$ has no indecomposable submodule with dimension
    vector $(2,1)$.

Let $B_i = B(e_i)$ and $B_{ij} = B(e_i+e_j)$ where $e_1,\dots,e_n$ is the canonical basis of $k^n$
and $i\neq j.$
The bristle $B_i$ is embedded as the submodule generated by $d(c_i)
= c_{ii}$ and the bristle $B_{ij}$ with $i\neq j$ is embedded as the submodule
generated by $d(c_i+c_j) = c_{ii}+c_{ij}+c_{jj}$. These are $\binom {n+1}2$
submodules and their sum is equal to $C$. In particular, we see that
$C$ is a bristled module. Since every element
of the socle of $C$ is obviously contained in a bristle, $C$ is a fully bristled
module. On the other hand, we also know that any bristle $B$ occurs as a unique
submodule of $C$, thus $\dim\Hom(B,C) = 1.$

It remains to be seen that $\End(C) = k.$ The assertion is clear for $n=1$,
thus we assume that $n\ge 2$. First, let us show that any non-zero endomorphism $f$
of $C$ is an automorphism. Namely, assume that $f$ is an endomorphism with non-zero kernel.
Let $U$ be a proper non-zero subspace of $C_2$ which belongs to the kernel of $f$.
Let $c\in C_2\setminus U.$ We claim that $C/U$ has an indecomposable  submodule $V$
with dimension vector $(2,1)$ whose socle is generated by $c+U$. Namely, if
$0 \neq u\in C$, then there are bristles $B$ and $B'$ with socle $kc$ and $k(c+u)$,
respectively, and the sum of the images of $B$ and $B'$ in $C/U$ is the requested
submodule. Now $f$ induces a homomorphism $C/U \to C$. Since
$\Hom(V,B) = 0$, it follows that also $c$ belongs to the kernel of $f$. Using
induction, we see that $f$ vanishes on $C_2$ and therefore $f = 0.$

Consider now an automorphism $f$ of $C$. Since $\dim\Hom(B,C) = 1$ for each bristle $B$,
we see that $f$ maps any submodule of $C$ which is a bristle into itself, thus
it maps any one-dimensional subspace of $C_2$ into itself. This shows that $f|C_2$
is a scalar multiplication. But this implies that $f$ itself is a scalar multiplication.
\hfill$\square$
	\bigskip
	{\bf Remark.}
	The map $d\:\Bbb PC_2 \to \Bbb P C_1$ is the degree 2 Veronese embedding.
	    \bigskip
{\bf (2.2)}
{\it There is a natural bijection between quadratic forms $\phi$ on $C_2$ and
homomorphisms $\phi'\:C \to S(1)$ of quiver representations,}
it is given as follows: We denote by $k[x_1,\dots, x_n]$ the polynomial ring
in $n$ variables,
and we evaluate these polynomials on $C_2$ via $x_r(\sum_i\lambda_i c_i)
 = \lambda_r)$; the space of quadratic forms on $C_2$ is 
the set of quadratic homogeneous polynomials, this space has the basis
$x_rx_s$ with $1 \le r \le s \le n$. We may evaluate this polynomial $x_rx_s$ 
on $C_1$ in the obvious way: $x_rx_s(\sum_{i\le j} \lambda_{ij}c_{ij}) = \lambda_{rs},$
this yields a linear transformation $\phi'\:C_1 \to k$, and we may interpret this
linear transformation as a homomorphism $C \to S(1)$ of quiver representations
which we also denote by $\phi'$.
           \medskip
The maximal submodules of $C$ are just the kernels of homomorphisms $C \to S(1)$.
Recall that a submodule $U$ of a module $M$ of finite length is said to be an {\it essential} 
submodule provided
that $U$ contains the socle of $M$. Since $C$ is a non-simple indecomposable 
$K(n)$-module, the socle of $C$ is equal to the radical of $C$, thus the essential
submodules of $C$ are just the intersections of maximal submodules, or, equivalently,
the intersections of the kernels of a set of homomorphisms $C \to S(1)$. The relationship
between the vanishing set for $\phi$ and the kernel of $\phi'$ is given as follows:
      \medskip
													 {\bf (2.3)} {\it Let $\phi$ be a quadratic form on $C_2$ and $c\in C_2$. Then $\phi(c) = 0$
if and only if $d(c)$ belongs to the kernel of $\phi'$.}
																		            \medskip
Proof: This is a tautological statement, since we have $\phi = \phi'\circ d$. Indeed, let us apply
$\phi$ and $\phi'\circ d$ to  $c = \sum_i \lambda_i c_i.$ If
$\phi$ is a basis element, say $\phi = x_rx_s$, then $x_rx_s(c) = \lambda_r\lambda_s$,
and $(\phi'\circ d)(c) = \phi'(d(c)) = x_rx_s(\sum_{i,j}\lambda_i\lambda_jc_{ij}) =
 \lambda_r\lambda_s.$ In general, $\phi$ is a linear combination of such basis elements,
and $\phi'$ is the corresponding linear combination.
			    	         \hfill$\square$
																											\medskip
{\bf (2.4)} The representation $C$ of $K(n)$ can be obtained also by looking
at the corresponding Beilinson algebra $A(n)$, this is the path algebra with quiver
$\Delta(n)$ with 3 vertices, say $1,2,3$, with $n$ arrows
$1\to 2$ labeled $\alpha_1,\dots,\alpha_n$ as well as $n$ arrows $2\to 3$, 
also labeled $\alpha_1,\dots,\alpha_n$, modulo the ``commutativity''
relations $\alpha_i\alpha_j = \alpha_j\alpha_i$ (whenever this
 makes sense).
$$
\hbox{\beginpicture
  \setcoordinatesystem units <3cm,2.5cm>
\put{$1$} at 0 0
\put{$2$} at 1 0
\put{$3$} at 2 0
\arr{0.2 0.1}{0.8 0.1}
\arr{0.2 -.1}{0.8 -.1}
\put{$\alpha_1$} at 0.5 0.2
\put{$\alpha_{n}$} at 0.5 -.2
\put{$\vdots$} at  0.5 0.035
\arr{1.2 0.1}{1.8 0.1}
\arr{1.2 -.1}{1.8 -.1}
\put{$\alpha_1$} at 1.5 0.2
\put{$\alpha_{n}$} at 1.5 -.2
\put{$\vdots$} at  1.5 0.035
\endpicture}
$$
For $i=1,2,3,$ let
$I_{A(n)}(i)$ be the injective envelope and $P_{A(n)}(i)$ the projective cover
of the $A(n)$-module $S(i)$.
If $Z$ is an $A(n)$-module, we denote by
$Z|[1,2]$ the restriction of $Z$ to the subquiver of type $K(n)$ with
vertices $1,2$. Then
$$
   C = I_{A(n)}(3)|[1,2].
   $$

Proof: Let $k[x_1,\dots,x_n]$ be the polynomial ring in $n$ variables $x_1,\dots,x_n$.
We may identify the $A(n)$-module $P_{A(n)}(1)$ with the subspace of
$k[x_1,\dots,x_n]$ generated by the monomials of degree $i-1$, for $1\le i \le 3$.
This shows that $P_{A(n)}(1)$ has dimension vector $(1,n,\binom{n+1}2)$.
Using duality, we see
that the dimension vector of $I_{A(n)}(3)$ is $(\binom{n+1}2,n,1)$.

On the other hand, we may extend $C$ to an $A(n)$-module by setting
$C_3 = k$, with maps
$\alpha_i:C_2 \to C_3$ defined by $\alpha_i(c_j) = 1$ if $i=j$ and zero otherwise
(obviously, the commutativity relations are satisfied). This extended module
has simple socle, namely the simple module $S(3)$,
thus it can be embedded into $I_{A(n)}(3)$, and the
dimension vectors show that the extended module is equal to $I_{A(n)}(3)$.
	    \bigskip\medskip
{\bf 3. Proof of Theorem 1$'$.}
		       \medskip
We assume now that $k$ is algebraically closed. We want to show that
any projective variety $\Cal V$ occurs as $\beta(M)$
for a reduced Kronecker module $M.$
It is well-known that we can realize $\Cal V$ as a closed subset of a
projective space, say $\Bbb P^{n-1}$, defined by a finite set of quadratic homogeneous
polynomials, say $q_1,\dots,q_m$. We consider the corresponding linear maps
$q'_i\:C \to S(1)$, where $C$ is the canonical $K(n)$-module and
put $M = \bigcap_i \Ker(q'_i)$.

Thus, let $\Cal V$ be the closed subset of the projective space $\Bbb P^{n-1}$,
defined by the vanishing of homogeneous polynomials $q_1,\dots,q_m$ of degree 2.
Let $\Lambda'$ be the factor
algebra of the Beilinson algebra $A(n)$ taking the elements $q_1,\dots,q_m$
as additional relations (obviously, these elements may be considered as
linear combinations of paths of length 2 in the quiver $\Delta(n)$). The
injective envelope $I_{\Lambda'}(3)$ of $S(3)$ as a $\Lambda'$-module
is just the submodule of $I_{A(n)}(3)$ which is obtained as the kernel
of the linear maps $q'_1,\dots,q'_m\:I_{A(n)}(3) \to S(1).$

Since $C$ is reduced, also its submodule $M$ is reduced. This completes the proof.
      \bigskip\bigskip
{\bf 4\. Open problems and remarks.}
      	    \medskip
{\bf (4.1) Problems.} Looking at Theorem 1 (and 1$'$),
a lot of questions may be asked. Let us mention just a few.
         \medskip
{\bf (a)} Theorem 1 asserts that any projective variety can be realized as
the eigenvector variety of some matrix pencil $P = (a,b;\alpha_1,\dots,\alpha_n)$,
or, equivalently, as the bristle variety of an $n$-Kronecker module.
The construction of $P$ given in the proof shows that even for quite innocent
projective varieties, the number $n$ may be large. Thus, we may ask {\it which
projective varieties occur as the eigenvector variety of a matrix pencil with
a fixed number $n$ of matrices,} or, equivalently, as the bristle variety of
an $n$-Kronecker module with $n$ fixed.

For $n = 1$, the reduced matrix pencils $P = (a,b;\alpha)$ are given by an arbitrary
injective linear transformation $\alpha\:k^a \to k^b$. It follows that any
non-zero element of $k^a$ is an eigenvector for $P$, thus
$\epsilon(P) = \Bbb P(k^a)$ is a projective space.

For $n = 2,$ one can use the well-known classification of the indecomposable
matrix pencils (see [4] or also [9]) in order to determine all
possible eigenvector varieties $\epsilon(P)$. 
Let us mention here which varieties occur as eigenvector varieties $\epsilon(P)$ where
$P$ is indecomposable (and different from $S(1)$).
The dimension vector of an indecomposable matrix pencil $P$ 
is of the form $(a,b)$ with $|a-b| \le 1$.
If $a < b$, then $\epsilon(P) = \emptyset.$ If $a = b$, then $\epsilon(P)$
consists of at most one point (if $k$ is algebraically closed, then $\epsilon(P)$
is non-empty). Finally, if $a > b$, then $\epsilon(P)$ is a non-singular quadratic
curve.

        Thus, the first case of real interest is the case $n = 3$: which
	projective varieties occur as the eigenvector variety of a pencil of
        $3$ matrices?
	          \medskip
	{\bf (b)} If we realize the projective variety $\Cal V$ as $\beta(M)$ for some
	reduced Kronecker module $M$, the elements of $\Cal V$ are the bristle
	submodules of $M$.
	Now certain elements of $\Cal V$ are of particular interest, for example the
        singular points. One may ask about specific properties of the
	        corresponding bristle submodules of $M$.
		       \medskip
{\bf (c)} Also, we may fix the dimension of $\Cal V$. Already the case of
     dimension $1$ seems to be of interest. Thus, one may ask:
{\it For which matrix pencils  $P$ is $\epsilon(P)$ a projective curve?}
     	  \bigskip
{\bf (4.2) Square matrices.} The main result of this note may be strengthened as
follows: {\it Let $k$ be an algebraically closed field. Any projective variety occurs
as the eigenvector variety of a reduced matrix pencil $P = (a,b;\alpha_1,\dots,\alpha_n)$
with square matrices $\alpha_i$.}
	       \medskip
Proof. Again, the proof will be given in terms of $n$-Kronecker modules.  
We will show that for any reduced
$n$-Kronecker module $M$ say with dimension vector $\bdim M = (a,b)$ 
there is a reduced $n$-Kronecker module $N$ with dimension vector
$\bdim N = (a',a')$ such that $\beta(M) = \beta(N)$ and $a' \le a$.
The proof relies on the following lemma.
	       	   \medskip
{\bf Lemma.} {\it Let $M$ be any $n$-Kronecker module, say with dimension vector $\bdim M = (a,b).$
If $M$ is the sum of its bristle submodules, then $M$ is the sum of $\ a\ $ bristle submodules.
As a consequence, we have $a\ge b$.}
		      \medskip
Proof. For the proof of the first assertion, we use induction on $a$. Let $M$ be
an $n$-Kronecker module with dimension vector $(a,b)$ and assume that $M$ is a sum of bristle
submodules. If $a = 0$, then also $b = 0$. Thus, assume that $a > 0.$
Since $M$ has finite length, $M$ is the sum of finitely many bristle
submodules, say of the bristle submodules $U_1,\dots, U_m$, and we assume that $m$ is minimal.
Of course, we must have $a \le m$. Let $M' = \sum_{i=1}^{m-1} U_i$ and
$\bdim M' = (a',b')$. Assume that
$a' = a$. Choose a vector space complement $V$ of $M'_2$ in $M_2$, thus $\dim V = b-b'$.
Then $M'' = (0,V)$ is a submodule of $M$ and $M = M'\oplus M''$, thus $M''$
is a direct summand of $M$. But
$\Hom(B,M'') = 0$ for any bristle. Since $M$ is generated by bristle submodules, $M''$
is generated by the images of maps $B \to M''$, where $B$ are bristles. It follows that
$M'' = 0$, thus $b' = b$, and therefore $M' = M.$ But this contradicts
the minimality of $m$. This shows that we must have $a' < a$.
By induction, $M'$ is the sum of $a'$
bristle submodules. Since $M$ is the sum of $M'$ and $U_m$, we see that $M$ is the sum
of $a'+1$ bristle submodules. The minimality of $m$ implies that $m \le a'+1$, thus
$m \le a'+1 \le a$. Altogether, we see that $m = a$, thus $M$ is the sum of $a$ bristle submodules.
This completes the proof of the first assertion.

In order to show the second assertion, we note that
the inclusion maps $U_i \to M$ yield a surjective map $\bigoplus_{i=1}^a U_i \to M$. Therefore
$$
 (a,a) = \bdim \bigoplus\nolimits_{i=1}^a U_i \ \ge\ \bdim M = (a,b).
$$
This shows that $a \ge b$. \hfill$\square$

 \medskip
Now, let $M$ be any reduced $n$-Kronecker module.
Let us denote by $M'$ the sum of all bristle submodules of $M$. If $\bdim M' = (a',b')$,
then $a' \le a$ (and $b' \le b$).
Since $M$ is reduced, also its submodule $M'$ is reduced, and, of course, we have
$$
  \beta(M') = \beta(M).
  $$
According to the Lemma, we have $a' \ge b'$. Let $N =
M' \oplus S(2)^{a'-b'}$. Since $M'$ is reduced, also $N$ is reduced.
Since $\Hom(B,S(2)) = 0$ for any bristle $B$, we see that any bristle submodule of $N$
is contained in $M'\oplus 0$, thus
$$
   \beta(N) = \beta(M').
$$
On the other hand, we have $\bdim N = (a',a')$.
   \hfill$\square$
	\bigskip
{\bf (4.3) Non-reduced matrix pencils.} It seems worthwhile to have a short look also at
non-reduced matrix pencils and non-reduced $n$-Kronecker modules. Any $n$-Kronecker module
$M$ can be written as a direct sum $M = M'\oplus M''$, where $M'$ is a direct sum
of copies of $S(1)$ and $M''$ is reduced; the submodule $M'$ is uniquely determined and,
of course, $M'' = M/M'$ is the maximal reduced factor module of $M$. In [11] we have
proposed to call $\eta(M) = \beta(M'')$ the {\it bristle variety of $M$}. In this way,
the bristle variety of any $n$-Kronecker module is defined and is a projective variety.

Looking at an arbitrary $n$-Kronecker module $M$, we also may consider the set
$\beta(M)$ of bristle submodules of $M$, this is an (open) subset of $\Bbb G_{(1,1)}(M)$.
As we have mentioned above, we have $\beta(M) = \Bbb G_{(1,1)}(M)$ in
case $M$ is reduced. This equality also holds true (for trivial reason) in case $M_2 = 0$,
since then both sets $\beta(M)$ and $\Bbb G_{(1,1)}(M)$ are empty.

Let us assume now that
$M$ is not reduced and that $M_2 \neq 0$. In this case, $M$ has a submodule isomorphic to $S(1)\oplus
S(2)$, thus $\beta(M)$ is a proper subset of $\Bbb G_{(1,1)}(M)$. In case
$\beta(M)$ is not empty, it is not closed in $\Bbb G_{(1,1)}(M)$ (thus not a projective
variety). Namely, let $U$ be a bristle submodule of $M$. Since we assume that $M$ is not
reduced, $M$ has a submodule $V$ isomorphic to $S(1)$ and, of course, $U\cap V = 0.$
The submodule $U\oplus V$ of $M$ shows that $\soc(U\oplus V)$
is a degeneration of $U$ (and $\soc(U\oplus V)$
belongs to $\Bbb G_{(1,1)}(M)$, but not to $\beta(M)$).
	\bigskip
	{\bf (4.4) Wild acyclic quivers.} The result presented in this paper can be used in order
	to show: {\it Let $k$ be an algebraically closed field. If $Q$ is any connected wild acyclic
	quiver with at least $3$ vertices, then any projective variety arises as a quiver Grassmannian
	of some representation of $Q$.} The proof will be given in [12]. It is an open question whether
	the condition that $Q$ has at least 3 vertices can be omitted (thus, whether any
	projective variety arises as a quiver Grassmannian of some $3$-Kronecker quiver).

   \bigskip\bigskip
{\bf 5\. Appendix. On the existence of eigenvectors.}
     \medskip
The aim of this appendix is to provide a translation of the main result of [11]
to the matrix pencil language, namely to reformulate this result as an existence
assertion for eigenvectors. We assume that $n\ge 3.$ 
	\medskip

{\bf (5.1)} We will say that a reduced matrix pencil $P = (a,b;\alpha_1,\dots,\alpha_n)$
has {\it sufficiently many eigenvectors} provided the eigenvectors of $P$ generate
the vector space $k^a$ (or, equivalently, provided the corresponding Kronecker module is
bristled). What we will see is that there is a huge class of matrix pencils
with sufficiently many eigenvectors.

If $v$ is an eigenvector of the reduced matrix pencil $P = (a,b;\alpha_1,\dots,\alpha_n)$,
and $w$ is a non-zero vector in
$\langle \alpha_1(v),\dots,\alpha_n(v)\rangle$, let $\alpha_i(v)
= \lambda_iw$, for $1\le i \le n$. The element
$\langle(\lambda_1,\dots,\lambda_n)\rangle\in E = \Bbb P(k^n)$ does not depend on
the choice of $w$ and will be called the {\it eigenvalue of $v$}. Of course,
equivalent eigenvectors have the same eigenvalue.

If $E'$ is a subset of $E$, we say that the reduced matrix pencil
$P = (a,b;\alpha_1,\dots,\alpha_n)$
has {\it sufficiently many $E'$-eigenvectors} provided
the eigenvectors with eigenvalue in $E'$ generate $k^a$. For $n\ge 3$,
we are going to exhibit a finite (!) subset $E_0$ of $E$ such that there
is a huge class of matrix pencils with sufficiently many $E_0$-eigenvectors.

As before, we will denote by $e_1,\dots,e_n$ the canonical basis of $k^n$.
If $1\le i,j\le n$ and $i\neq j$, let $e_{i,j} = e_i+e_j$. We denote by
$E_0$ the subset of $E$ given by the following elements: first, 
$\langle e_{n-1}\rangle$ and $\langle e_n\rangle$, second $\langle e_{i,i+1}\rangle$
for $1\le i < n$, and finally $\langle e_{1,n}\rangle.$
Note that $E_0$ is a set of cardinality $n+2$. 
     \medskip
Let us now define the reflection functor $\sigma$ as introduced by Bernstein-Gelfand-Ponomarev.
We recall that matrix pencils $P = (a,b;\alpha_1,\dots,\alpha_n)$ and
$P' = (a',b';\alpha'_1,\dots,\alpha'_n)$ are said to be {\it equivalent} provided
$a = a', b= b'$ and there are invertible linear transformations $\beta\:k^a \to k^a$
and $\gamma\:k^b \to k^b$ such that $\gamma\alpha_i\beta = \alpha'_i$ for all
$1\le i \le n$. The reflection functor $\sigma$ sends the (equivalence class of the) 
matrix pencil $P = (a,b;\alpha_1,\dots,\alpha_n)$
to the following equivalence class $\sigma P$ of matrix pencils: Let $U$ be the set
of elements $u = (u_1,\dots,u_{na})$ in $k^{na}$ such that
$\sum_{i=1}^n \alpha_i(u_{(i-1)a+1},\dots,u_{(i-1)a+a}) = 0.$ Let $z$ be the dimension of $U$.
Write the inclusion map $U \to k^{na} = k^a\times \cdots\times k^a$
in terms of some basis of $U$ and the
canonical basis of $k^a\times \cdots \times k^a$;
we obtain in this way a sequence of $n$ $(a\times z)$-matrices $\omega_1,\dots,\omega_n.$
Then $\sigma P$ is the equivalence class of $(z,a;\omega_1,\dots,\omega_n).$
     	    \medskip
{\bf Theorem 2.} {\it Let $n\ge 3.$ If $P = (a,b;\alpha_1,\dots,\alpha_n)$
is any matrix pencil, then, for $t\gg 1$, the matrix pencil
$\sigma^tP$ has sufficiently many $E_0$-eigenvectors.}
		    	            \medskip
Let us stress that Theorem 2 is essentially a
reformulation of part of the main theorem in [11]. 
We note the
following: if $B(\lambda)$ is a bristle $n$-Kronecker module, then any non-zero
vector in $B(\lambda)_1$ is an eigenvector with eigenvalue $\lambda$, and,
conversely, any eigenvector of a matrix pencil $P$ with eigenvalue $\lambda$
gives rise to an embedding of $B(\lambda)$ into $P$, considered as a Kronecker
module. In particular, the set $\Cal B_0$ of bristles which is
used in [11] corresponds bijectively to the set $E_0 \subset \Bbb P(k^n)$.
In the paper [11], instead of the functor $\sigma$ the Auslander-Reiten translation $\tau$
of the category of $n$-Kronecker modules was used. But according to Gabriel [3], we have
$\tau = \sigma^2.$ Thus it is sufficient to to apply (1.3)(b) of [11] both for $M(P)$ and 
$\sigma M(P)$.
																										    \medskip
													
{\bf Remark.} The reader should be aware that for $n = 2$, the corresponding assertion
is not valid, typical examples are the matrix pairs
$\left(2,2;\left[\smallmatrix 1 & 0 \cr
                              0 & 1 \endsmallmatrix \right],
      			       \left[\smallmatrix \lambda & 1 \cr
	                      0 & \lambda \endsmallmatrix \right]\right)$
			       with $\lambda\in k.$
	    \bigskip
{\bf (5.2)} Let us mention under what condition one may recover $P$ from
 $\sigma^tP$, where $t$ is a natural number. Following [9],
 we say that an indecomposable matrix pencil $P$ is {\it preprojective} provided
 $\sigma^tP = 0$ for some natural number $t$. If $P = (a,b;\alpha_1,\dots,\alpha_n)$
 is preprojective, then $a < b$ and $a^2+b^2-nab = 1$ and any indecomposable matrix pencil
 with this dimension vector $(a,b)$ is equivalent to $P$
(also, any pair $(a,b)$ of natural
numbers with $a < b$ and $a^2+b^2-nab = 1$ arises in this way). It follows that for
$n\ge 2$, there are countably  many equivalence classes of preprojective matrix pencils,
and it is easy to see that for $n\ge 2$,
{\it the preprojective matrix pencils have no eigenvectors.}

If $P$ is an indecomposable matrix pencil which is not preprojective, then
we can recover $P$ from
$\sigma^t P$ by using a corresponding functor $\sigma^{-t}$ (the left adjoint for
$\sigma^t$): namely, $P$ and $\sigma^{-t}\sigma^t P$ are equivalent matrix pencils.
Since Theorem 2 asserts that
$\sigma^tP$ has sufficiently many eigenvectors for $t \gg 1$,
we see that {\it any indecomposable
matrix pencil which is not preprojective is of the form $\sigma^{-t}P'$
for some $t$, where $P'$ is a matrix pencil with sufficiently many eigenvectors.}
  \bigskip\bigskip

{\bf Acknowledgment.}
The author is grateful to Birge Huisgen-Zimmermann and Lutz Hille for useful
comments concerning the final presentation.
  \bigskip\bigskip

{\bf 6\. References.}
       \medskip
 \item{[1]} Bongartz, K., Huisgen-Zimmermann, B.:
      Varieties of uniserial representations IV.
     Kinship to geometric quotients
     Trans Amer. Math. Soc. 353 (2001), 2091--2113.
\item{[2]} Derksen, H., Huisgen-Zimmermann, B., Weymann, J.:
        Top stable degenerations of finite
  dimensional representations II.  Advances in Mathematics. 259 (2014), 730--765
\item{[3]} Gabriel, P.: Auslander-Reiten sequences and representation-finite algebras.
             In: {\it Representation Theory I.} Spinger Lecture Notes in Math\. 831 (1980), 
             1--71.
\item{[4]} Gantmacher, F. R.: Matrizentheorie, Springer Verlag (1986).
\item{[5]} Hille, L.: Tilting line bundles and moduli of thin sincere representations
   	    of quivers. An. St. Univ. Ovidius Constantza 4 (1996), 76--82.
\item{[6]} Hille, L.: Moduli of representations, quiver Grassmannians and Hilbert schemes.
	   	       arXiv: 1505.06008.
\item{[7]} Huisgen-Zimmermann, B.: The geometry of uniserial representations of
	   		              finite dimensional algebras I. J. 
           Pure Appl. Algebra 127 (1998),     39--72.
\item{[8]} Reineke, M.: Every projective variety is a quiver Grassmannian.
	   Algebra and Represent. Theory 16 (2013), no. 5, 1313--1314.
\item{[9]} Ringel, C\. M.: {\it Tame algebras and integral quadratic forms.}
	   	            Springer Lecture Notes in Math\. 1099 (1984).
\item{[10]} Ringel, C\. M.: Quiver Grassmannians and Auslander varieties
	   	               for wild algebras. J.Algebra 402 (2014), 351--357.
\item{[11]} Ringel, C\. M\.: Kronecker modules generated by modules of length 2.
   To appear in {\it Representations of Algebras.} Contemporary Mathematics. Amer.~Math.~Soc.,
      	     	    arXiv:1612.07679.
\item{[12]} Ringel, C\. M\.: Quiver Grassmannians for wild acyclic quivers. arXiv:1703.08782.
		    			    	 \bigskip\bigskip
{\rmk

Claus Michael Ringel \medskip

Fakult\"at f\"ur Mathematik, Universit\"at Bielefeld\par
D-33501 Bielefeld, Germany
\medskip

Department of Mathematics, Shanghai Jiao Tong University \par
Shanghai 200240, P. R. China.
	  \medskip

e-mail: \ttk ringel\@math.uni-bielefeld.de \par}

\bye